
\magnification=\magstep1
\def\forces{\parallel\!\!\! -}


\def\hexnumber#1{\ifcase#1 0\or1\or2\or3\or4\or5\or6\or7\or8\or9\or
        A\or B\or C\or D\or E\or F\fi }

\font\teneuf=eufm10
\font\seveneuf=eufm7
\font\fiveeuf=eufm5
\newfam\euffam
\textfont\euffam=\teneuf
\scriptfont\euffam=\seveneuf
\scriptscriptfont\euffam=\fiveeuf


\font\tenmsx=msam10
\font\sevenmsx=msam7
\font\fivemsx=msam5
\font\tenmsy=msbm10
\font\sevenmsy=msbm7
\font\fivemsy=msbm5
\newfam\msxfam
\newfam\msyfam
\textfont\msxfam=\tenmsx  \scriptfont\msxfam=\sevenmsx
  \scriptscriptfont\msxfam=\fivemsx
\textfont\msyfam=\tenmsy  \scriptfont\msyfam=\sevenmsy
  \scriptscriptfont\msyfam=\fivemsy
\edef\msx{\hexnumber\msxfam}

\mathchardef\upharpoonright="0\msx16

\def\qed{{\vcenter{\hrule height.4pt \hbox{\vrule width.4pt height5pt
 \kern5pt \vrule width.4pt} \hrule height.4pt}}}
\def\notin{{\in}\kern-5.5pt / \kern1pt}
\def\ok{\vbox{\hrule height 8pt width 8pt depth -7.4pt
    \hbox{\vrule width 0.6pt height 7.4pt \kern 7.4pt \vrule width 0.6pt height 7.4pt}
    \hrule height 0.6pt width 8pt}}
\def\nt{{\leq}\kern-1.5pt \vrule height 6.5pt width.8pt depth-0.5pt \kern 1pt}
\def\sd{{\times}\kern-2pt \vrule height 5pt width.6pt depth0pt \kern1pt}
\def\zp#1{{\hochss Y}\kern-3pt$_{#1}$\kern-1pt}

\def\extend{ \hat {\; \; } }
\def\sm{{\smallskip}}

\def\no{{\noindent}}

\def\sub{\subseteq}

\def \o {\omega }

\font\capit=cmcsc10 scaled\magstep0

\overfullrule=0pt
\openup1.5\jot
\def\sm{\smallskip}

\def \o {\omega } \def \sub {\subseteq } \def \s {\sigma }

\def \r {\rho } \def \k {\kappa } 
\def\g{\gamma}
\def\l{\lambda } \def\th {\theta} 
\def\d{\delta}

\font\sgross=cmbx10 scaled \magstep2
  \def\z{\zeta}
\def\a{\alpha } 
\def\b{\beta}
\def\m{\mu}
\def \e {\varepsilon}
\def\epsilon {\varepsilon}
\def\l {\lambda}
\def\bigprod{\prod}

\noindent {\sgross On incomparability and related cardinal functions
on ultraproducts of Boolean algebras}\footnote{}{{1991 {\sl
Mathematics Subject Classification.} Primary 06E05, 03E35; Secondary
03E55. Key words: Boolean algebra, ultraproduct, cardinal
characteristic, large cardinal}}

\bigskip 

\centerline{{\capit Saharon Shelah}\footnote{$^1$}{Supported by grant
NSF-DMS97-04477 publication 677.} {\capit and Otmar
Spinas}\footnote{$^2$}{Supported by grant
2124-045702.95/1 of the Swiss National Science Foundation and by grant
NSF-DMS97-04477.} }

\bigskip 

{\narrower

{\noindent ABSTRACT: Let $C$ denote any of the following cardinal
characteristics of Boolean algebras: incomparability, spread,
character, $\pi$-character, hereditary Lindel\"of number, hereditary
density. It is shown to be consistent that
there exists a sequence $\langle B_i:i<\kappa\rangle$ of Boolean
algebras and an ultrafilter $D$ on $\k$ such that
$$C(\prod_{i<\k}B_i/D)<|\prod_{i<\k}C(B_i)/D|.$$ This answers a number
of problems posed in [M].

}}

\bigskip

{\bf Introduction}

\bigskip

For a number of cardinal characteristics $C$ of Boolean algebras it
makes sense to ask whether it is consistent to have a sequence
$\langle B_i:i<\kappa\rangle$ of Boolean algebras and an ultrafilter
$D$ on $\k$ such that
$$C(\prod_{i<\k}B_i/D)<|\prod_{i<\k}C(B_i)/D|.$$ For $C$ being the
length this was proved in [MSh]. The same method of proof can be used
to get the analogous thing for $C$ being any one of the following:
incomparability ($Inc$), spread ($s$), character ($\chi$),
$\pi$-character ($\pi\chi$), hereditary Lindel\"of number ($hL$),
hereditary density ($hd$). This ansers problems 47, 48, 52, 56, 60 of
[M]. For irredundancy (Monk's problem 25) this will be done in a
subsequent paper of the first author. We won't define these notions
here, as they are very clearly defined on pp. 2,3 in [M]. We assume
that the reader has a good knowledge of [MSh] and [Mg].

For $C$ a cardinal function of Boolean algebras which is defined as
the supremum of all cardinals which have a certain property, we define
$C^+$ as the least cardinal $\k$ such that this property fails for
every cardinal $\lambda\geq \k$. Note that all cardinal functions
mentioned above are of this form. For $\chi$, note that in [M]
$\chi(B)$ has been defined as the minimal $\k$ such that every
ultrafilter on $B$ can be generated by $\k$ elements. Clearly,
$\chi(B)$ can be equivalently defined as $$\sup\{\chi(U):U\hbox{ is an
ultrafilter on }B\},$$ where $\chi(U)$ is the minimal size of a
generating subset of $U$. 

In $\S 1$ below we deal with incomparability. In $\S 2$ we shall show
that the results for all the other characteristics can be deduced from
this relatively easily.

A key notion for the proofs is that of $\mu$-entangled linear order,
$\mu$ being a cardinal (see definition 1.1 below). The reason for this
is the following observation of Shelah (see [M, p.225]), where
$Int(I)$ denotes the intervall algebra of some linear order $I$,
i.e. the subalgebra of ${\cal P}(I)$ generated by the half-open
intervals of the form $[a,b)$, $a,b\in I$.

\smallskip

{\bf Fact.} {\it Let $\mu$ be a regular uncountable cardinal and let
$I$ be a linear order. The following are equivalent:

\item{(1)} $I$ is $\mu$-entangled.

\item{(2)} There is no incomparable subset of $Int(I)$ of size $\mu$.}

\smallskip

For $b$ a member of some Boolean algebra $B$, $b^1$ denotes $b$ and $b^0$
denotes the complement of $b$. By $Ult(B)$ we denote the Stone space
of $B$.

\bigskip

{\bf 1. Incomparability}

\bigskip

{\bf Definition 1.1.} Let $(I,<)$ be a linear order and let $D\subseteq
{\cal P}(\kappa)$ for some infinite cardinal $\kappa$.

(1) $(I,<)$ is called $(\delta,\gamma)$-entangled, where $\gamma
<\delta$ are cardinals, if for every family $\langle
t_{\alpha,\epsilon }: \alpha <\delta ,\epsilon
<\epsilon(\ast)<\gamma\rangle$ of pairwise distinct members of $I$ and
for every $u\subseteq \epsilon(\ast)$ there exist $\alpha <\beta
<\delta$ such that $$\forall \epsilon <\epsilon(\ast)\quad t_{\alpha
,\epsilon}<t_{\beta,\epsilon} \Leftrightarrow \epsilon\in u.$$ If
$\gamma =\o$ we say that $(I,<)$ is $\delta$-entangled.

(2) $(I,<)$ is called $(\delta ,D)$-entangled if for every sequence
$\langle t_{\alpha,\epsilon,l }: \alpha <\delta ,\epsilon \in A,l<n
\rangle$ of pairwise distinct members of $I$, where $A\in D$ and
$n<\o$, and for every $u\subseteq n$ there exist $\alpha
<\beta <\delta$ and $B\subseteq A$, $B\in D$, such that
$$\forall \epsilon \in B\forall l<n \quad t_{\alpha
,\epsilon,l}<t_{\beta,\epsilon,l} \Leftrightarrow l\in u.$$

\smallskip

Note that if $\kappa <\gamma$ and $(I,<)$ is
$(\delta,\gamma)$-entangled then $(I,<)$ is $(\delta ,D)$-entangled
for every $D\sub {\cal P}(\kappa)$.

In the sequel, if $\gamma<\delta$ are regular cardinals, by
$C(\gamma,\delta)$ we denote the partial order to add $\delta $ Cohen
subsets to $\gamma$. More precisely, 
$$C(\gamma,\delta)=\bigprod_{i<\delta } Q_i \qquad (<\gamma
\hbox{-support product})$$ where $Q_i=({^{<\gamma}2}, \supseteq)$. Clearly
$Q$ is $\gamma$-directedly-closed.

\smallskip

{\bf Lemma 1.2.}  {\it Let $\gamma<\delta$ be regular cardinals such
that $\forall \a<\d\quad \a^{<\gamma}<\delta$. Let
$L=\{\eta_i:i<\delta\} $ be $C(\gamma,\delta)$-generic. Letting
$<_{lex}$ denote the lexicographic order on $^\gamma 2$, we have that
$(L,<_{lex})$ is a $(\delta ,\gamma)$-entangled linear order.}

\sm

{\it Proof:} Suppose $p\forces_{C(\gamma ,\delta)} $ ``$\langle \dot
\tau(\a,\epsilon): \a<\d ,\epsilon<\epsilon(\ast)<\gamma \rangle$ is a
sequence of pairwise distinct ordinals below $\delta$ such that the
family $\langle \dot \eta_{\dot \tau(\a,\epsilon)}: \a<\d
,\epsilon<\epsilon(\ast)<\gamma \rangle$ contradicts
$(\delta,\gamma)$-entangledness of $L$, witnessed by set
$u\sub\epsilon(\ast)$''. Here $\dot \eta _i$ is the canonical name for
the $i$th Cohen subset of $\gamma$.

As $C(\gamma ,\delta)$ does not add new ordinal sequences of length
$<\gamma$ we may assume that $u\in V$. For the same reason, for each
$\a<\delta$ we may pick $p_\a\leq p$ such that $p_\a$ decides the
value of $\langle \dot \tau
(\a,\epsilon):\epsilon<\epsilon(\ast)\rangle$, say as $\langle \tau
(\a,\epsilon):\epsilon<\epsilon(\ast)\rangle$. By the $\Delta$-system
Lemma and some thinning out there exists $Y\in [\delta]^\delta$ and
$r\in [\delta]^{<\gamma}$ such that for all $\a,\b\in Y$, $\a\ne\b$,
we have:

\sm

\item{(i)} dom$(p_\a)$ $\cap$ dom$(p_\b) =r$,

\item{(ii)} $\{\tau (\a,\epsilon):\epsilon<\epsilon(\ast)\}\sub$
dom$(p_\a)$,

\item{(iii)} $\langle p_\a(i):\a\in Y\rangle$ is constant for every
$i\in r$,


\item{(iv)} $\langle p(\tau(\a,\epsilon)):\a\in Y\rangle$ is constant
for every $\epsilon<\e(\ast)$.
\sm

Pick $\a,\b\in Y$ with $\a <\b$. Note that $p_\a, p_\b$ are
compatible, and hence $\{\tau
(\a,\epsilon):\epsilon<\epsilon(\ast)\}\cap \{\tau
(\b,\epsilon):\epsilon<\epsilon(\ast)\}=\emptyset$.

Define $q\leq p_\a,p_\b$ by dom$(q)=$ dom$(p_\a)$ $\cup$ dom$(p_\b)$
and $$q(i)=\cases{p_\a(i)\extend 0 &$i\in \{\tau
(\a,\epsilon):\epsilon\in u\}$ \cr 
p_\a(i)\extend 1 &$i\in \{\tau
(\a,\epsilon):\epsilon\in \e(\ast)\setminus u\}$ \cr
p_\b(i)\extend 0 &$i\in \{\tau
(\b,\epsilon):\epsilon\in \e(\ast)\setminus u\}$ \cr
p_\b(i)\extend 1 &$i\in \{\tau
(\b,\epsilon):\epsilon\in u\}$ \cr 
p_\a(i) &$i\in$ dom$(p)\setminus \{\tau(\a,\e):\e<\e(\ast)\}$\cr
p_\b(i) &otherwise.\cr}$$ Then clearly $$q\forces\; (\forall
\e<\e(\ast)) \eta_{\dot \tau(\alpha ,\e)}<_{lex} \eta_{\dot \tau(\b
,\e)}\Leftrightarrow \e\in u,$$ a contradiction. \hfill $\qed$

\bigskip

Assume GCH. Let $\mu$ be a supercompact cardinal and let $\k <\m$ be a
measurable cardinal. Fix $D$ a normal measure on $\k$. By [L] we may
assume that the supercompactness of $\m$ cannot be destroyed by any
$\m$-directedly-closed forcing.

For any ordinal $\alpha$ let $F(\a)$ denote the least inaccessible
cardinal above $\alpha$, if it exists. We assume that $F(\m)$ exists
and denote it with $\lambda$.

Let $Q=C(\mu ,\lambda)$. Clearly $Q$ is $\m$-directedly-closed
and $V^Q\models 2^\m=\l$.

Work in $V^Q$. Let $U$ be a normal fine measure on $[H((2^\l)^+)]^{<\m}$. By
Lemma 1.2, $H((2^\l)^+)\models $ ``there exists a $(\l,\m)$-entangled
linear order on $\l$''. Therefore the set $A$ of all $a\in 
[H((2^\l)^+)]^{<\m}$ such that 

\sm

\item{(1)} $(a,\in)\prec (H((2^\l)^+),\in)$,

\item{(2)} $a\cap \m$ is measurable,

\item{(3)} the Mostowski collapse of $a$ is $H((2^{F(a\cap \m)})^+)$,

\item{(4)} $H((2^{F(a\cap \m)})^+)\models$ there is a $(F(a\cap\m)
,a\cap \m)$-entangled linear order on $F(a\cap \m)$, call it $J^*_{a\cap\m}$,

\item{(5)} $H((2^{F(a\cap \m)})^+) \models 2^{a\cap\m}=F(a\cap \m)$

\sm \no belongs to $U$.

By well-known arguments on large cardinals and elementary embeddings
we can build a sequence $\bar U=\langle U_\a:\a<\k\rangle$ of normal
measures on $\m$ such that

\sm

\item{(a)} $\a <\b\Rightarrow U_\a <U_\b $ (i.e. $U_\a\in $
Ult$(V,U_\b)$),

\item{(b)} $\{ a\cap\m:a\in A\}\in U_\a$ for all $\a <\k$.

\no The main fact which is used for this is the following lemma which
goes back to [SRK]. We thank James Cummings for reconstructing the
proof for us. 

\smallskip

{\bf Lemma 1.3} {\it For all $a\in V_{\mu +2}$ there exists a normal
measure $U$ on $\mu$ such that $a\in \, Ult(V,U)$ and $\{ a\cap\m:a\in
A\}\in U $.} 

\sm

{\it Sketch of proof:} Let $j:V\rightarrow M$ be the elementary
embedding defined by the normal fine measure $U$ above. Fix $<$ a
wellordering of $V_\mu$ and let $$<^*=j(<)\upharpoonright
V_{\lambda^+}.$$ Assuming that the Lemma is false, let $b\in
V_{\mu+2}$ be the $<^*$-minimal counterexample. Let $\bar U=\{ B\sub
\mu:\mu\in j(B)\}$ be the normal measure on $\mu$ induced by $j$ and
let $i:V\rightarrow N$ be the corresponding elementary embedding. As
usual we have another elementary embedding $k:N\rightarrow M$, defined
by $k([f]_{\bar U})=j(f)(\mu)$, such that $j=k\circ i$ (see [J,
p.312]). As $V_{\mu+2}\sub M$ we have that $$M\models b\hbox{ is the
}j(<)\hbox{-minimal counterexample.}$$ By elementarity there must
exist $\bar b\in N$ such that $k(\bar b)=b$ and $$N\models \bar
b\hbox{ is the }i(<)\hbox{-minimal counterexample.}$$ Note that
$b\not\in N$, as $b$ is a counterexample. Also note that $V_{\mu
+1}\sub \; ran(k)$. As $k$ is simply the inverse of the transitive
collapse map on $ran(k)$, we conclude $\bar b =b$ and hence $b\in N$,
a contradiction. \hfill $\qed$

\sm

Let $Q(\bar U)$ denote Magidor's forcing to change the cofinality of
$\m$ to $\k$ by adding a normal sequence $\langle \m_i:i<\k\rangle$
cofinal in $\m$. Fix such a $Q(\bar U)$-generic sequence with
$\m_0>2^\k$. We let $$\m_i'=\m_{\o i},\, \th _i=F(\m_{i+1}),\,
\l_i=F(\m_{\o i}), \, J_i=J^*_{\m_{i+1}}.$$ 

\sm

{\bf Lemma 1.4} {\it For every $i<\kappa$, $V^{Q\ast Q(\bar U)}
\models \, ``F(\mu_{i+1})=\theta_i $ and $J_i$ is
$(\theta_i,D)$-entangled''.} 

\sm

{\it Proof:} Work in $V^Q$. Let $\langle \dot \mu _i:i<\kappa\rangle $
be a $Q(\bar U)$-name for the generic sequence. Fix $i<\kappa$. Let
$p\in Q(\bar U)$ such that $p$ decides $\dot \mu _j$ as $\mu _j$ for
$j\in \{ i,i+1,i+2\}$. We may assume that the domain of the first
coordinate of $p$ is $\{ i,i+1,i+2\}$. By the main arguments of [Mg],
especially [Mg, Lemma 5.3], it follows that forcing $Q(\bar U)$ below
$p$ factors as $P^i_{\mu _i}\ast Q^i$, where $P^i_{\mu_i}$ is the union
of $\mu_i$ many $\mu_i$-directed suborders each of them of size $\leq
2^{\mu_i}$, and $Q^i$ does not add new subsets to $\mu_{i+2}$. Hence
clearly $V^{Q\ast Q(\bar U)} \models \, F(\mu_{i+1})=\theta_i $.

Now suppose $p\forces \, ``\langle \dot t_{\alpha,\epsilon,l} :\alpha
<\theta_i,\epsilon \in A,l<n\rangle $ is a one-to-one family of
elements of $J_i$''. By [Mg, Lemma 4.6], for each $\alpha <\theta_i$ we
can find $p_\alpha \leq p$ such that $p_\alpha $ and $p$ have the same
first coordinate and for all $\epsilon\in A$ and $l<n$ there exists
$w_{\alpha ,\epsilon,l}\in [i]^{<\o}$ such that below $p_\alpha $, the
value of $\dot t_{\alpha,\epsilon,l}$ depends only on the value of
$\langle \dot \mu_j:j\in w_{\alpha,\epsilon,l}\rangle $. As $D$ is
$\kappa$-complete and $i<\kappa$, there exists $B_{\alpha,l}\in D$ and
$w_{\alpha,l}$ such that $w_{\alpha,\epsilon,l}=w_{\alpha,l}$ for all
$\epsilon\in B_{\alpha,l}$ and $l<n$. Let $w_\alpha
=\bigcup_{l<n}w_{\alpha,l}$, $B_\alpha =\bigcap_{l<n}B_{\alpha,l}$. As
$2^\kappa <\theta_i$ and $i<\theta_i$ we can find $Y\sub \theta_i$ of
size $\theta_i$, $w^*\in [i]^{<\o}$ and $B\in D$ such that $B_\alpha
=B$ and $w^*=w_\alpha $ for all $\alpha\in Y$. Let $v$ be the domain
of the first coordinate of any $p_\alpha$. By [Mg, Lemma 3.3], for each
$\alpha\in Y$ we can find $p_\alpha'\leq p_\alpha$ such that the
domain of the first coordinate of $p_\alpha '$ is $w^*\cup v$. Then
$p_\alpha'$ decides $\langle \dot \mu_j:j\in w^*\rangle $, say as
$\langle \mu^\alpha_j:j\in w^*\rangle $, and hence $p_\alpha'$ decides
$\langle \dot t_{\alpha,\epsilon,l}:\epsilon\in B,l<n\rangle$, say as
$\langle t_{\alpha,\epsilon,l}:\epsilon\in B,l<n\rangle$. Note that
this sequence is one-to-one. As $\mu_i<\theta_i$ we can find $Y'\sub
Y$ of size $\theta_i$ and $\langle \mu_j:j\in w^*\rangle $ such that
$\langle \mu^\alpha_j:j\in w^*\rangle = \langle \mu_j:j\in
w^*\rangle$ and $(p_\alpha )_{i+2}=(p_\beta)_{i+2}$ (we use the
notation of [Mg, p.67]), for all $\alpha ,\beta \in Y'$. By [Mg, Lemma
4.1] it follows that $p_\alpha$ and $p_\beta$ are compatible for all
$\alpha,\beta\in Y'$. It follows that $\langle  t_{\alpha,\epsilon,l} :\alpha
\in Y',\epsilon \in B,l<n\rangle $ is a one-to-one family. Applying
$(\theta_i,D)$-entangledness of $J_i$ in $V^Q$, for any $u\sub n$ we
obtain $B'\in D$, $B'\sub B$ and $\alpha <\beta$, $\alpha,\beta\in
Y'$, such that for all $\epsilon\in B'$, for all $l<n$,
$t_{\alpha,\epsilon,l}<t_{\beta,\epsilon,l}\Leftrightarrow l\in
u$. What we have shown suffices to prove the Lemma. \hfill $\qed$

\sm

For every $i<\kappa$ we define a linear order $I_i\sub \bigprod_{j<\o
i}\theta_j$ as follows: For every $i'<\o i$ fix a family $\langle
A^\r:\r\in \bigprod_{j<i'}\theta_j\rangle $ of pairwise disjoint
subsets of $\theta_{i'}\cap \; Card$, each of them of cardinality
$\theta_{i'}$. This is possible as $|\bigprod_{j<i'}\theta_j|<\theta_{i'}$
and $\theta_{i'}$ is a regular limit cardinal. Let $I_i$ be the set of
all $\eta\in\bigprod_{j<\o i}\theta_j$ such that for all $j<\omega i$,
$\eta(j)\in A^{\eta\upharpoonright j}$. Define a linear order $<_i$ on
$I_i$ as follows: For distinct $\eta,\nu\in I_i$ let $\e=\min\{ j<\o
i: \eta(j)\ne \nu(j)\}$. Now let $$\eta<_i\nu \Leftrightarrow
\cases{\e \hbox{ is even and }\eta(\e) <_{J_\e} \nu(\e), \hbox{ or}&
\cr \e \hbox{ is odd and }\eta(\e) < \nu(\e).& \cr}$$

\sm

We claim that in $I_i$ we can choose a one-to-one family $\langle
\eta^i_{\z ,\e}:\e\leq\z<\l_i \rangle$ such that the following hold:

\sm

\itemitem{(a)}
$\forall\z_1<\z_2\forall\e_1\leq\z_1\forall\e_2\leq\z_2\quad
\eta^i_{\z_1 ,\e_1}<_{J^{bd}_{\o i}} \eta^i_{\z_2 ,\e_2},$

\sm

\itemitem{(b)} $\langle \eta^i_{\z ,0}:\z<\l_i\rangle$ is cofinal in
$\bigprod_{j<\o i}\theta_j /{J^{bd}_{\o i}}$,

\sm

\itemitem{(c)} the mapping $\langle (\eta^i_{\z,2\e},\eta^i_{\z,2\e
+1}):\e <\z\rangle $ is $<_i$-preserving.

\sm

Here ${J^{bd}_{\o i}}$ denotes the ideal of bounded subsets of $\o i$.
For the construction of such a family remember from [MSh] that $
\bigprod_{j<\o i}\theta_j /{J^{bd}_{\o i}} $ has true cofinality
$\lambda _i$. Clearly, in $I_i$ we can find a family $\langle
\eta^i_\z :\z<\l_i\rangle$ which is increasing and cofinal in
$\bigprod_{j<\o i}\theta_j /{J^{bd}_{\o i}}$ and satisfies
$\eta^i_\zeta(j)\cdot 3 <\eta^i_{\zeta +1}(j)$ for almost all $j<\o
i$. Now let $\z<\l_i$ and $2\e<\z$. Define $\eta_{\z,2\e}^i$ and
$\eta_{\z,2\e+1}^i$ by letting
$$\eta^i_{\z,2\e +l}(j) =\cases{\eta_\z(j) & if $j$ is even, \cr
\eta_\z(j) + \eta_\epsilon (j) & if $j$ is odd and $l=0$, \cr
\eta_\zeta(j) + \eta_\e(j)\cdot 2 & if $j$ is odd and $l=1$. \cr }$$
It is easy to see that this definition works.

\sm

{\bf Lemma 1.5.} {\it In $V^{Q\ast Q(\bar U)}$ the following holds:
Whenever $\langle J_i:i<\kappa\rangle$ is a family such that for every
$i$, $J_i$ is a $(\theta_i,D)$-entangled linear order on $\theta _i$
and $I_i$ is defined as above, then $(I_i,<_i)$ is
$(\l_i,D)$-entangled but not $(\l',\aleph_0)$-entangled for any
$\l'<\l_i$.}

\sm

{\it Proof:} The last statement easily follows from the existence of
the family $\langle \eta^i_{\z ,\e}:\e\leq\z<\l_i \rangle$. Let
$\langle t_{\a,\e,l}: \a<\l_i,\e\in A ,l<n\rangle $ be a family of
pairwise distinct members of $I_i$, where $A\in D$ and $n<\o$. Hence
$$t_{\a,\e,l}=\eta^i_{\z(\a,\e,l) ,\nu(\a,\e,l)}$$ for some
$\nu(\a,\e,l)\leq \z(\a,\e,l)<\l_i$. Fix $\a<\lambda_i$. As $i<\k$
there is $A_\a'\in D$ and $i^*_\a <\o i$ such that for all distinct
$\e,\e'\in A'_\a$ and $l,m<n$ do we have $t_{\a,\e,l}\upharpoonright
i^*_\a \ne t_{\a,\e ',l}\upharpoonright i^*_\a$,
$t_{\a,\e,l}\upharpoonright i^*_\a \ne t_{\a,\e ,m}\upharpoonright
i^*_\a$ and $t_{\a,\e,l}\upharpoonright i^*_\a \ne t_{\a,\e
',m}\upharpoonright i^*_\a$. As $2^\k <\l_i$ we may assume that
$\langle A'_\a:\a<\l_i\rangle $ and $\langle i^*_\a: \a<\l_i\rangle $
are constant, say with values $A^*$, $i^*$. As in [Sh462, Claim 3.1.1]
one shows that there must exist cofinally many even $j\in (i^*,\o i)$
such that for every $\xi <\theta _j$ there is $\a <\l_i$ with the
property $\forall\e\in A^*\forall l<n\quad t_{\a,\e,l}(j)>\xi$. Fix
such $j$. Construct an increasing sequence $\langle \a(\nu):\nu
<\theta_j\rangle$ such that $$\forall\nu<\rho<\theta_j\forall
\e,\e'\in A'\forall l,m<n\quad t_{\a(\nu),\e,l}(j) < t_{\a(\rho ),\e
',m}(j).$$ As $(\bigprod_{l<j}\theta_l)^\k <\theta_j$, we may assume
that the sequence $\langle \langle t_{\a(\nu),\e,l}\upharpoonright j:
\e\in A^*,l<n\rangle : \nu <\theta_j\rangle $ is constant. Note that by
construction, $$\langle t_{\a(\nu),\e,l}(j): \nu <\theta_j,\e\in
A^*,l<n\rangle $$ is a sequence of pairwise distinct members. We can
apply $(\theta_j,D)$-entangledness of $J_j$ and, for given $u\subseteq
n$, we get $B\in D$, $B\sub A^*$, and $\nu <\r<\theta_j$ such that
$$\forall\e\in B\forall l<n\quad t_{\a(\nu),\e,l}(j) <_{J_i}
t_{\a(\r),\e,l}(j)\Leftrightarrow l\in u.$$ By construction we
conclude that $$\forall\e\in B\forall l<n\quad t_{\a(\nu),\e,l} <_i
t_{\a(\r),\e,l}\Leftrightarrow l\in u.$$ \hfill $\qed$

\sm

{\bf Lemma 1.6.} {\sl Letting $I=\bigprod_{i<\k}I_i/D$, $I$ is $(\l
,\aleph_0)$-entangled in $V^{Q\ast Q(\bar U)\ast \, Coll(\mu^+,<\l)}$.}

\sm

{\it Proof:} By Lemmas 1.4 and 1.5 and as $Coll(\mu^+,<\l)$ does not
add new subsets to $\mu$, in $V^{Q\ast Q(\bar U)\ast \,
Coll(\mu^+,<\l)}$ it is true that $I_i$ is $(\lambda_i,D)$-entangled
and $D$ is a normal fine measure on $\kappa $. Moreover, note that
$\bigprod_{i<\k}\l_i/D$ has order-type $\lambda$ in $V^{Q\ast Q(\bar
U)\ast \, Coll(\mu^+,<\l)}$. This is true because it holds in
$V^{Q\ast Q(\bar U) }$ by [MSh] and because $Coll(\mu^+,<\l)$ does
not add new functions to $\bigprod_{i<\k}\l_i$.  As $V^{Q\ast Q(\bar
U)\ast \, Coll(\mu^+,<\l)}\models \lambda =\mu ^{++}$ we have that the
cofinality of $\bigprod_{i<\k}\l_i/D$ is $\lambda$.

Let $\langle t^l_\a:\a<\l,l<n\rangle $, $n<\o$, be a
family of pairwise distinct elements of $I$. So $t^l_\a$ is of the
form $$t^l_\a =\langle
\eta^i_{\z_i(\a,l),\varepsilon_i(\a,l)}:i<\k\rangle /D,$$ where
$\eta^i_{\z_i(\a,l),\varepsilon_i(\a,l)}\in I_i$.  By the above
observations, wlog we may assume that 

\sm

\itemitem{$(\ast_1)$} $\langle \langle \z_i(\a,l):i<\k\rangle/D
:\a<\l\rangle $ is increasing and cofinal in $\bigprod_{i<\k}\l_i/D$,
for every $l<n$.

\sm

For every $\a<\l$ and $i<\k$ there is $j<\o i$ such that for every
$l<m<n$, if $\langle \z_i(\a,l),\varepsilon_i(\a,l)\rangle \ne
\langle \z_i(\a,m),\varepsilon_i(\a,m)\rangle $ then
$$\eta^i_{\z_i(\a,l),\varepsilon_i(\a,l)} \upharpoonright j\ne
\eta^i_{\z_i(\a,m),\varepsilon_i(\a,m)} \upharpoonright j.$$ By Los'
Theorem and since $D$ is normal, there exist $B_\a'\in D$ and
$j_\a'<\k$ such that for all $i\in B'_\a$ and $l<m<n$ we have
$$\eta^i_{\z_i(\a,l),\varepsilon_i(\a,l)} \upharpoonright j'_\a \ne
\eta^i_{\z_i(\a,m),\varepsilon_i(\a,m)} \upharpoonright j'_\a.$$ As
$2^\k<\l$, wlog we may assume that
\sm

\itemitem{$(\ast_2)$} there are $B^1\in D$ and $j'<\k$ such that for
all $\a<\l$, $i\in B^1$ and $l<m<n$
$$\eta^i_{\z_i(\a,l),\varepsilon_i(\a,l)} \upharpoonright j' \ne
\eta^i_{\z_i(\a,m),\varepsilon_i(\a,m)} \upharpoonright j'.$$
\sm

Moreover we have
\sm

\itemitem{$(\ast_3)$} there exist $B^2\in D$, $B^2\sub B^1$, and
$\langle j^2_i:i\in B^2\rangle$ such that $j^2_i<\o i$ and for every
$g\in \bigprod_{i\in B^2}\l_i$, $f_i\in \bigprod\{\theta_j:j_i^2\leq
j<\o i\}$ and $\a<\l$ we can find $\b\in (\a,\l)$ such that for every
$i\in B^2$, $j^2_i\leq j<\o i$ and $l<n$ we have $g(i)<\z_i(\b,l)$ and
$$f_i(j)<\eta ^i_{\z_i(\b,l),\varepsilon_i(\b,l)}(j).$$

If $(\ast_3)$
failed, for every candidate $y=\langle B^y,\langle j_i^y:i\in
B^y\rangle\rangle $ to satisfy $(\ast_3)$ we had $g^y$, $\langle
f^y_i:i\in B^y\rangle$, $\a^y$ which witness that $y$ does not satisfy
$(\ast_3)$. Note that there are only $2^\k$ candidates. Let $$\a=\sup\{
\a^y:y\hbox{ is a candidate}\}$$ and $$f_i(j)=\sup\{f^y_i(j):y\hbox{
is a candidate and }j\in\, dom(f^y_i)\}.$$ As there are only $2^\k$
candidates we have $\a<\l$ and $f_i(j)<\theta_j$. We can choose
$\b_i<\l_i$ such that $\b_i>\z_i(\a,l)$ for every $l<n$ and
$$f_i<_{J^{bd}_{\o i}}\eta^i_{\b_i,\varepsilon}$$ for every
$\varepsilon \leq \b_i$.
Finally we define $g\in\bigprod_{i<\k}\l_i$ by letting
$$g(i)=\sup\{g^y(i): y\hbox{ is a candidate and }i\in B\}\cup\{\b_i+1\}.$$

By $(\ast_1)$ we can find $\gamma\in (\a,\l)$ and $B\in D$ such that
$g(i)<\langle \z_i(\gamma,l):i<\k\rangle$ for all $i\in B$ and
$l<n$. By construction, for every $i\in B$ there is $j_i<\o i$ such
that for all $j_i\leq j<\o i$ and $l<n$
$$f_j(j)<\eta^i_{\z_i(\g,l)\e_i(\gamma,l) }(j).$$ Then $y=\langle
B,\langle j_i:i\in B\rangle\rangle $ is a candidate which contradicts
the definition of $\a,$ $\langle f_i:i<\k\rangle$, $g$. This finishes
the proof ok $(\ast_3)$.

As $D$ is normal, wlog we may assume that in $(\ast)_3$, $\langle
j^2_i:i\in B^3\rangle $ is constant with value $j^2<\k$. Now choose
$i^*\in B^3$ even with $\max\{ j^1,j^2\} <i^*$. Using $(\ast _3)$ it
is straightforward to find an increasing sequence $\langle
\a(\nu):\nu <\theta_{i^*}\rangle$ in $\lambda$ such that for all $i\in
B^3\setminus i^*+1$ and $l,m<n$ we have
$$\eta^i_{\zeta_i(\a(\nu),l)\e_i(\a(\nu),l)} (i^*)<
\eta^i_{\zeta_i(\a(\nu +1),l)\e_i(\a(\nu+1),l)} .$$ As
$(\prod_{j<i^*}\theta_j)^\k < \theta_{i^*}$, wlog we may assume that
$$\langle \langle
\eta^i_{\zeta_i(\a(\nu),l)\e_i(\a(\nu),l)}\upharpoonright i^*: i\in
B^3\setminus (i^*+1),l<n\rangle: \nu <\theta_{i^*}\rangle $$ is
constant. By construction we have that, letting
$$s_{\nu,i,l}=\eta^i_{\zeta_i(\a(\nu),l)\e_i(\a(\nu),l)} (i^*), $$ 
$$\langle s_{\nu,i,l}: \nu <\theta_{i^*},i\in
B^3\setminus(i^*+1),l<n\rangle $$ is a sequence of pairwise distinct
members of $I_{i^*}$. Hence by Lemma 1.5, for every $u\subseteq n$ we can
find $\nu<\xi<\theta_{i^*}$ and $A\in D$, $A\subseteq B^3\setminus
(i^*+1)$ such that for all $i\in A$ and $l<n$ we have
$$s_{\nu,i,l}<s_{\xi,i,l}\Leftrightarrow l\in u.$$ This implies
$$t^l_{\a(\nu)}<t^l_{\a(\xi)} \Leftrightarrow l\in u,$$ which finishes
the proof. \hfill $\qed$

\sm

As a corollary we obtain the following:

\sm

{\bf Theorem 1.7} {\it For $i<\kappa$ let $I_i$ be the linear order
defined above and let $B_i=$ $Int(I_i)$. In the model $V^{Q\ast Q(\bar
U)\ast \, Coll(\mu^+,<\l)}$ the following hold:

\sm

(i) $Inc(B_i) =\, Inc^+(B_i)=\lambda _i$ for all
$i<\kappa$, and hence $\bigprod_{i<\kappa}Inc(B_i)/D =\lambda
=\mu^{++}$,

\sm

(ii) $Inc^+(\bigprod_{i<\kappa} B_i/D) \leq\lambda $ and
hence $Inc(\bigprod_{i<\kappa} B_i/D) \leq \mu^{+}.$}

\sm

{\it Proof:} (i) follows from the fact mentioned in the introduction
and Lemma 1.5. Note that Lemma 1.5 holds also in $V^{Q\ast Q(\bar
U)\ast \, Coll(\mu^+,<\l)}$ as $Coll(\mu^+,<\l)$ does not add new
subset of $\mu$. 

(ii) follows from the same fact, by Lemma 1.6 and by the fact that
$\bigprod_{i<\kappa}B_i/D$ is isomorphic to Int$\,\bigprod_{i<\kappa}I_i/D$.
This last fact holds by Los' Theorem and as $D$ is
$\aleph_1$-complete. \hfill $\qed$

\bigskip

{\bf 2. Other characteristics}

\bigskip

{\bf Definition 2.1.} If $(I,<)$ is a linear order, by $Sq(I)$ we
denote the Boolean subalgebra of $({\cal P}(I^2),\sub)$ generated by
sets of the form $$X_{a,b}=\{ (a',b')\in I^2: a'<a\hbox{ and
}b'<b\},$$ for $a,b\in I$.

\sm

Recall that a sequence $\langle y_\alpha:\alpha <\lambda\rangle $ of
elements of some Boolean algebra is {\it left-separated} iff for every
$\alpha <\lambda$, $y_\alpha$ does not belong to $Id\langle
y_\b:\b>\a\rangle$, the ideal generated by $\langle
y_\b:\b>\a\rangle$. Similary, $\langle y_\alpha:\alpha <\lambda\rangle
$ is {\it right-separated} if for every $\alpha <\lambda$, $y_\alpha$
does not belong to $Id\langle y_\b:\b<\a\rangle$, the ideal generated
by $\langle y_\b:\b<\a\rangle$.

\sm

{\bf Lemma 2.2.} {\it Suppose $(I,<)$ is a $\l$-entangled linear order,
where $\l>\o$ is regular. Then $Sq(I)$ has neither a left-separated
nor a right-separated sequence of length $\l$. }

\sm

{\it Proof:} We prove the Lemma only for right-separated sequences. The
proof for left-separated sequences is similar. Suppose $\langle
y_\a:\a<\l\rangle$ is a right-separated sequence in Sq$(I)$. We shall
obtain a contradiction. Each $y_\a$ is a finite union of finite
intersections of sets of the form $X_{a,b}$ or $-X_{a,b}$. One of
these finite intersections does not belong to $Id\langle
y_\b:\b<\a\rangle$. Hence wlog we may assume that each $y_\a$ is such
a finite intersection. As cf$(\l)>\o$, wlog there exist $n<\o$ and
$\eta:n\rightarrow 2$ such that $$y_\a=
\bigcap_{l<n}X^{\eta(l)}_{a(\a,l),\b(\a,l)} $$ for some
$a(\a,l),\b(\a,l)\in I$, for all $\a<\l$.

\sm

{\it Case I:} $\exists l<n\quad \eta(l)=1$.

\sm

As the intersection of any two sets of the form $X_{a,b}$ has the same
form, wlog we may assume that $\eta(0)=1$ and $\eta(l)=0$ for all
$0<l<n$. We may also assume that $0<l<l'<n$ implies $a(\a,l)\ne
a(\a,l')$, $b(\a,l)\ne b(\a,l')$ and $a(\a,l)< a(\a,l')\Leftrightarrow
b(\a,l)> b(\a,l')$, for all $\a <\lambda$. Otherwise we could choose a
smaller $n$. Hence we have two subcases according to whether
$a(\a,1)<\dots <a(\a,n-1)$ and $b(\a,1)>\dots >b(\a,n-1)$ or
$a(\a,1)>\dots >a(\a,n-1)$ and $b(\a,1)<\dots <b(\a,n-1)$ holds. We
assume the first alternative holds. The second one is symmetric.

For fixed $\a<\l$ define the following sets:

\sm

\itemitem{$z_0=$} $X_{a(\a,0),b(\a,0)}-X_{a(\a,n-1),b(\a,0)},$

\sm

\itemitem{$z_1=$} $X_{a(\a,n-1),b(\a,0)}-X_{a(\a,n-1),b(\a,n-1)}
-X_{a(\a,n-2),b(\a,0)} ,$

\sm

\itemitem{$\dots$}

\sm

\itemitem{$z_{n-2}=$} $X_{a(\a,2),b(\a,0)}-X_{a(\a,2),b(\a,2)}
-X_{a(\a,1),b(\a,0)} ,$

\sm

\itemitem{$z_{n-1}=$} $X_{a(\a,1),b(\a,0)}-X_{a(\a,1),b(\a,1)}.$

\sm

Note that $y_\a=\bigcup _{j<n}z_j$. Hence there exists $j<n$ such that
$z_j\not\in \;Id\langle y_\b:\b<\a\rangle$. Wlog we may assume that
$j$ is the same for all $\a <\l$ and that $y_\a=z_j$ for all $\a
<\l$. Then $y_\a$ has the form $X_{a,b}-X_{a',b}-X_{a,b'}$ or
$X_{a,b}-X_{a',b}$ or $X_{a,b}-X_{a,b'}$ for some $a'<a$ and
$b'<b$. Let us assume $y_\a$ is of the first form. The others are even
easier to handle. Hence we have $$y_\a
=X_{c(\a,0),d(\a,0)}-X_{c(\a,0),d(\a,1)}-X_{c(\a,1),d(\a,0)},$$ where
$c(\a,1)< c(\a,0)$ and $d(\a,1)< d(\a,0)$.

\sm

Choose $F\sub 2\times 2$ maximal such that there exist
$\sigma:F\rightarrow I$ and cofinally many $\a\in\l$ with the property
that $(0,j)\in F$ implies $c(\a,j)=\sigma(0,j)$ and $(1,j)\in F$
implies $d(\a,j)=\sigma(1,j)$ for all $j<2$. Wlog we may assume that
the above holds for all $\a <\l$ and that for all $\a<\b<\l$ and
$(i,j)\in 2\times 2\setminus F$, if $i=0$ then $c(\a,j)\ne c(\b,j)$
and if $i=1$ then $d(\a,j)\ne d(\b,j)$. Depending on $F$ we have 16
cases to consider. However we consider only the case $F=\emptyset$, as
the others are similar.

We have more subcases to consider according to the order-type of the
sequence $\langle c(\a,0),c(\a,1),d(\a,0),d(\a,1)\rangle $. Wlog we
may assume that it does not depend on $\a$. We only work through two
typical examples. Let us first assume that this sequence consists of
pairwise distinct elements. As we assumed $F=\emptyset$ we conclude
that $\langle c(\a,j),d(\a,j):\a<\l,j<2\rangle$ is a family of
pairwise distinct elements. By $\l$-entangledness of $I$ we get
$\a>\b$ such that $c(\a,0)<c(\b,0)$, $d(\a,0)<d(\b,0)$,
$c(\a,1)>c(\b,1)$ and $d(\a,1)>d(\b,1)$. We conclude $y_\a\leq y_\b$,
a contradiction. Now suppose $c(\a,0)=d(\a,0)<c(\a,1)<d(\a,1)$. In
this case the family $\langle c(\a,j),d(\a,1):\a<\l,j<2\rangle$
consists of pairwise distinct elements. By $\l$-entangledness we
obtain $\a>\b$ such that $c(\a,0)<c(\b,0)$, $c(\a,1)>c(\b,1)$ and
$d(\a,1)>d(\b,1)$. Again we conclude $y_\a\leq y_\b$, a
contradiction. The other cases are similar. 

\sm

{\it Case II:} $\forall l<n\quad \eta(l)=0$.

\sm

Again we may assume that $a(\a,0)<a(\a,1)<\dots <a(\a,n-1)$ and
$b(\a,0)>b(\a,1)>\dots >b(\a,n-1)$ for all $\a <\lambda$. Notice that
wlog we may assume that $$X^0_{a(\a,n-1),b(\a,0)}\not\in \; Id\langle
y_\beta:\b<\a\rangle $$ for all $\a <\lambda$, as otherwise we may
replace $y_\a$ by $y_\a\cap X_{a(\a,n-1),b(\a,0)}$ and proceed as in
Case I. Hence wlog $$y_\a = X^0_{a(\a,n-1),b(\a,0)}$$ for all $\a
<\lambda$. Let $a_\a =a(\a,n-1), b_\a =b(\a,0)$. Clearly, if $a_\a
=a_\b$ for some $\a <\b$ then $b_\a <b_\b$, as otherwise $y_\a \leq
y_\b$. Similarly, $b_\a =b_\b$ implies $a_\a <a_\b$. As a
$\lambda$-entangled linear order does not have any increasing or
decreasing sequences of length $\lambda$, wlog we may assume that both
families $\langle a_\a:\a <\lambda\rangle$ and $\langle b_\a:\a
<\lambda\rangle$ are one-to-one. By a similar argument we may assume
that $a_\a\ne b_\a$ for all $\a <\l$ and also that $a_\a\ne
b_\b$ for all $\a\ne\b$. We can apply $\lambda$-entangledness of $I$ to the
family $\langle a_\a,b_\a:\a<\lambda\rangle$ and get some $\a >\b$
such that $a_\a>a_\b$ and $b_\a>b_\b$. Hence $y_\a\leq y_\b$, a
contradiction. \hfill $\qed$

\sm

{\bf Lemma 2.3.} {\it Let $(I,<)$ be a linear order and $\mu$ a
cardinal such that there exist $\{ (a_\a,b_\a):\a<\m\}\sub I^2$ and
$c\in I$ with the property that $a_\a\ne a_\b$, $b_\a<c$ and that
$a_\a<a_\b$ implies $b_\a<b_\b$ for all $\a,\b<\mu$, $\a\ne\b$. Then
$s^+(Sq(I))> \mu$ holds. }

\sm

{\it Proof:} Let $y_\a =X_{a_\a,c}-X_{a_\a,b_\a}$, for $\a<\mu$. Note
that $$y_\a\not\leq \bigcup_{\b\in F}y_\b$$ for all $\a<\m$ and finite
$F\sub\mu$ with $\a\not\in F$. Indeed, let $F_0=\{\b\in
F:a_\a<a_\b\}$, $F_1=F\setminus F_0$, let $\b_0$ be the subscript of
the smallest $a_\b$, $\b\in F_0$ and let $\b_1$ be the subscript of the
largest $a_\b$, $\b\in F_1$. Then $y_\a\setminus\bigcup_{\b\in F}y_\b=
y_\a\setminus (y_{\b_0}\cup y_{\b_1})$. As $(a_\a,b_\a)\in
y_\a\setminus (y_{\b_0}\cup y_{\b_1})$ we are done. Hence there exists
a family of ultrafilters $\langle U_\a:\a<\m\rangle$ with $y_\a\in
U_\a$ and $-y_\b\in U_\a$ for all $\a\ne\b$. Then $\langle
U_\a:\a<\m\rangle$ is a discrete set of cardinality $\m$ in the Stone
space of $Sq(I)$. \hfill $\qed$

\sm

{\bf Corollary 2.4.} {\it Using the notation of $\S 1$, letting
$B_i=\, Sq(I_i)$ for $i<\k$, in the model $V^{Q\ast Q(\bar
U)\ast \, Coll(\mu^+,<\l)}$ the following hold:

\sm

(i) $s( B_i) =\, s^+(B_i)=\, hL(B_i)=\, hL^+(B_i)=\,
hd(B_i)=\, hd^+(B_i)=\lambda _i$ for all $i<\kappa$, and hence
$\bigprod_{i<\kappa}s(B_i)/D =\bigprod_{i<\kappa}hL(B_i)/D
=\bigprod_{i<\kappa}hd(B_i)/D =\lambda =\mu^{++}$,

\sm

(ii) $hL^+ (\bigprod_{i<\kappa} B_i/D) =\, hd^+
(\bigprod_{i<\kappa} B_i/D) \leq\lambda $ and hence
$s(\bigprod_{i<\kappa} B_i/D)$, $hL (\bigprod_{i<\kappa} B_i/D)$ and
$hd(\bigprod_{i<\kappa} B_i/D)$ are all at most $\mu^{+}.$} \hfill $\qed$

\sm

{\it Proof:} We first prove (i). The proofs of Theorem 6.7 and Lemma 6.8
in [M] show that for every Boolean algebra $B$, if $hd(B)=\kappa$,
$\kappa$ being regular and infinite, then $hd(B)$ is attained
(i.e. there exists a subspace $X\sub\; Ult(B)$ with $d(B)=\kappa$) iff
$B$ has a left-separated sequence of length $\kappa$. Similarly, the
proof of Theorem 15.1 in [M] shows that if $hL(B)$ is regular and infinite,
then $hL(B)=\kappa$ is attained iff $B$ has a right-separated sequence
of length $\kappa$. As trivially $s^+(B)\leq \min\{ hL^+(B),\;
hd^+(B)\}$ and hence $s(B)\leq \min\{ hL(B),\; hd(B)\}$ holds, we
conclude that all cardinal coefficients of $B_i$ mentioned in (i) are
at most $\lambda _i$. That they are at least $\lambda _i$ follows from
Lemma 2.3, the construction of $I_i$ and the trivial fact that every
linear order of cardinality $\mu^+$, for some cardinal $\mu$, has a
subset of size $\mu$ which has an upper bound. 

In order to prove (ii) note that by Los' Theorem and
$\aleph_1$-completeness of $D$ we have that $\bigprod_{i<\k}B_i/D$ is
isomorphic to $Sq(\bigprod_{i<\k}I_i/D)$. By Lemmas 1.6 and 2.2 and
the previous argument we get (ii). \hfill $\qed$

\sm

{\bf Definition 2.5.} Let $\langle y_\a:\a <\lambda\rangle$ be a
one-to-one enumeration of some infinite linear order $(J,<_J)$. Define
a linear order $(L(J),<)$ by letting $L(J)=\{ (y_\a,\b):\a <\l ,\,
\b<\a\}$ and $$(y_\a,\b)<(y_\a',\b')\Leftrightarrow \, (y_\a<_J
y_{\a'})\vee (y_\a =y_{\a'}\wedge \b<\b').$$

\sm

{\bf Lemma 2.6.} {\it Let $\s$ be an infinite, regular cardinal which
is not the successor of a singular cardinal. Let $(J,<_J)$ be a
linear order of size $\lambda$ which does not have any increasing or
decreasing chain of length $\l$. Then $$\chi^+(Int(L(J)))
=\pi\chi^+(Int(L(J)))=\l$$ holds.}

\sm

{\it Proof:} As trivially $\pi\chi^+(B)\leq\chi^+(B)$ holds for every
Boolean algebra $B$, it suffices to show $\chi^+(Int(L(J)))\leq\l$ and
$\pi\chi^+(Int(L(J)))\geq\l$. Let $U$ be an ultrafilter on
$Int(L(J))$. Let $$L_U=\{ z\in L(J):(-\infty ,z)\in U\}.$$ Clearly
$L_U$ is a (possibly empty) end-segment of $L(J)$. It is
straightforward to see that $$\chi(U)\leq \,cf(L\setminus L_U) +\,
cf(L_U^*),$$ where the cofinality of a linear order is the minimal
length of a well-ordered cofinal subset, and $L_U^*$ is the inverse
order of $L_U$. We claim that $cf(L\setminus L_U) +\,
cf(L_U^*)<\l$. Let us first consider $cf(L\setminus L_U)$. If
$(L\setminus L_U)\cap J\times\{0\}$ is unbounded in $L\setminus L_U$
then $cf(L\setminus L_U)$ equals the cofinality of some well-ordered
increasing chain in $J$, which is assumed to be $<\l$. Otherwise
$L\setminus L_U \sub \{ (y_\b,\gamma):\b\leq\a ,\,\gamma <\b\}$ for
some $\a<\l$. Then $cf(L\setminus L_U)\leq |\a|<\l$. We conclude
$\chi^+(Int(L(J)))\leq\l$. 

In order to prove $\pi\chi^+(Int(L(J)))\geq\l$ let $\sigma <\l$ be
regular. Let $U$ be the ultrafilter on $Int(L(J))$ generated by the
intervals $$[(y_{\sigma +1},\a), (y_{\sigma+1 },\sigma)), \quad
\a<\sigma .$$ Now let $Y\subseteq \, Int(L(J))\setminus \{ 0\}$ be
dense in $U$. If $|U|<\sigma$ there exists $y\in Y$ such that $y\sub 
[(y_{\sigma +1},\a), (y_{\sigma+1 },\sigma))$ holds for cofinally many
$\a <\sigma$. This is clearly impossible. \hfill $\qed$

\smallskip

{\bf Corollary 2.7.} {\it Using the notation of $\S 1$ and definition
2.5, letting $B_i= \;Int(L(I_i))$ for $i<\kappa$, in the model
$V^{Q\ast Q(\bar U)\ast \, Coll(\mu^+,<\l)}$ the following hold:

\sm

(i) $\pi( B_i) =\, \pi^+(B_i)=\, \pi\chi(B_i)=\,
\pi\chi^+(B_i)=\lambda _i$ for all $i<\kappa$, and hence
$\bigprod_{i<\kappa}\chi(B_i)/D =\bigprod_{i<\kappa}\chi^+(B_i)/D
=\lambda =\mu^{++}$,

\sm

(ii) $\chi^+ (\bigprod_{i<\kappa} B_i/D) =\, \pi\chi^+
(\bigprod_{i<\kappa} B_i/D) =\lambda $ and hence
$\chi(\bigprod_{i<\kappa} B_i/D)=\;\pi\chi (\bigprod_{i<\kappa}
B_i/D)= \mu^{+}.$} \hfill $\qed$

\sm

{\it Proof:} First note that $I_i$, $i<\kappa$, has a dense subset of
size $\mu_{\o i}$. Indeed, for each $s\in \bigcup_{j'<\o
i}\bigprod_{j<j'}\theta _j$ choose $\eta_s\in I_i$ with $s\sub
\eta_s$ if this is possible. It is easy to see that the collection of
all these $\eta_s$ is dense in $I_i$. As there are only $\m_{\o i}$
many $s$ we are done. Hence clearly $I_i$ does not have a well-ordered
increasing or decreasing chain of length $\lambda _i$. Hence by Lemma 2.6
we have (i). By Los' Theorem and $\aleph_1$-completeness of $D$ we
have that $\bigprod_{i<\kappa} B_i/D$ is isomorphic to
$Int(L(\bigprod_{i<\k}I_i/D ))$. 
By Los'
Theorem again and as $\bigprod_{i<\k}\lambda_i =\lambda$, it follows
that $\bigprod_{i<\k}I_i/D$ does not have a well-ordered increasing or
decreasing chain of length $\lambda$. By Lemma 2.6 we conclude
(ii). \hfill $\qed$

\bigskip

\centerline{{\bf References}}

\smallskip

\item{[J]}{{\capit T. Jech,} {\bf Set theory,} Academic Press, New
York (1978).}

\item{[L]} {\capit R. Laver,} {\it Making the supercompactness of
$\kappa $ indestructible under $\kappa $-directed closed forcing},
{\bf Israel J. Math. 29} (1978), no. 4, 385--388.

\item{[Mg]} {\capit M. Magidor,} {\it Changing cofinality of
cardinals}, {\bf Fund. Math. 99} (1978), 61-71.

\item{[MSh]} {\capit M. Magidor and S. Shelah,} {\it Length of Boolean
algebras and ultraproducts}, {\bf Math. Japonica}, to appear.

\item{[M]} {\capit J.D. Monk,} {\bf Cardinal Invariants on Boolean
Algebras}, Birkh\"auser Verlag, Basel, 1996.

\item{[SRK]} {\capit R. Solovay, W. Reinhardt and A. Kanamori,} {\it
Strong axioms of infinity and elementary embeddings} {\bf
Ann. Math. Logic 13} (1978), no. 1, 73--116.

\item{[Sh462]} {\capit S. Shelah,} {\bf Cardinal arithmetic}, Oxford
Logic Guides, 29. Oxford Science Publications. The Clarendon Press,
Oxford University Press, New York, 1994.

\bigskip

\smallskip

{\it Addresses: First author:}
\smallskip

\item{}{Institute of Mathematics, Hebrew University, Givat Ram, 91904
Jerusalem, ISRAEL}

\item{}Department of Mathematics, Rutgers University, New Brunswick,
NJ 08903, USA

{\it e-mail:} shelah@math.huji.ac.il

\smallskip

{\it Second author:}
\smallskip

Mathematik, ETH-Zentrum HG G33.3, 8092 Z\"urich, SWITZERLAND

{\it Phone:} +41 1 632 34 05

{\it Fax:} +41 1 632 10 85

{\it e-mail:} spinas@math.ethz.ch

\bye